\documentclass[12pt]{article}
\usepackage[centertags]{amsmath}
\usepackage{amsfonts}
\usepackage{amssymb}
\usepackage{graphicx}
\usepackage{amssymb}
\usepackage[verbose,colorlinks=true,naturalnames=true,linkcolor=blue,]
{hyperref}

\newcommand{\sign}{\mathop{\rm sign}\nolimits}
\begin{document}

\title{Multiparameter Quantum Deformations of Jordanian Type for
Lie Superalgebras
\footnote{This work is supported by the grants RFBR-05-01-01086 and ANR
grant NT05-241455GIPM. 
}}

\author{V.N. Tolstoy
\\
Institute of Nuclear Physics, 
Moscow State University,\\
119992 Moscow, Russia\footnote{\uppercase{E}mail:tolstoy@nucl-th.sinp.msu.ru}
}


\maketitle

\begin{abstract}{ We discuss quantum deformations of Jordanian type for Lie
superalgebras. These deformations are described by twisting functions with support from
Borel subalgebras and they are multiparameter in the general case. The total twists are
presented  in explicit form for the Lie superalgebras $\mathfrak{sl}(m|n)$ and
$\mathfrak{osp}(1|2n)$. We show also that the classical $r$-matrix for a light-cone
deformation of $D=4$ super-Poincare algebra is of Jordanian type and a corresponding
twist is given in explicit form. }
\end{abstract}

\section{Introduction}
The Drinfeld's quantum group theory roughly includes two classes
of Hopf algebras: quasitriangular and triangular. The (standard)
$q$-deformation of simple Lie algebras belongs to the first class.
The simplest example of the triangular (non-standard) deformation
is the Jordanian deformation of  $\mathfrak{sl}(2)$.
In the case of simple Lie algebras of rank $\geq2$ some
non-standard deformations were constructed by  Kulish, Lyakhovsky
et al.\cite{KLM}--\cite{AKL}.
These deformations are described by twisting functions (which are
extensions of the Jordanian twist) with support from Borel
subalgebras, and they are multiparameter in the general case.
We call their as the deformations of Jordanian type.
Total twists of Jordanian type were constructed for all complex Lie
algebras of the classical series $A_n$, $B_n$, $C_n$ and $D_n$.

In this paper we discuss quantum deformations of Jordanian type for
Lie superalgebras. The total twists are presented  in explicit form for the
Lie superalgebras $\mathfrak{sl}(m|n)$ and $\mathfrak{osp}(1|2n)$.
We show also that the classical $r$-matrix for a light-cone
deformation of $D=4$ super-Poincare algebra is of Jordanian type
and a corresponding twist is given in explicit form.

\section{Classical $r$-matrices of Jordanian type}

Let $\mathfrak{g}$ be any finite-dimensional complex
simple Lie superalgebra  then
$\mathfrak{g}=\mathfrak{n_-^{}}\oplus\mathfrak{h}\oplus\mathfrak{n_+^{}}$,
where $\mathfrak{n_{\pm}^{}}$ are maximal nilpotent subalgebras
and $\mathfrak{h}$ is a Cartan subalgebra. The subalgebra
$\mathfrak{n_+^{}}$ ($\mathfrak{n_-^{}}$) is generated by the
positive (negative) root vectors $e_{\beta}^{}$ ($e_{-\beta}^{}$ )
for all $\beta\in\triangle_+(\mathfrak{g})$. The symbol
$\mathfrak{b}_{+}^{}$ will denote the Borel subalgebra of
$\mathfrak{g}$,
$\mathfrak{b}_{+}^{}:=\mathfrak{h}\oplus\mathfrak{n}_+^{}$.
Let $\theta$ be a maximal root of $\mathfrak{g}$, and let
a Cartan element $h_{\theta}^{}\in \mathfrak{h}$ and a root vector
$e_{\theta}^{}\in \mathfrak{n}_{+}$ satisfies the relation
\begin{equation} \label{jr1}
[h_{\theta},\,e_{\theta}]=e_{\theta}.
\end{equation}
The elements $h_{\theta}^{}$ and  $e_{\theta}^{}$ are homogeneous,
i.e.
\begin{equation} \label{jr2}
\deg(h_{\theta})=0,\;\; \deg(e_{\theta})=0,\,{\rm or}\;\,1~.
\end{equation}
Moreover, let homogeneous elements $e_{\pm i}$ indexed by the symbols
$i$ and $-i$, ($i=1,2,\ldots,N$), satisfy the relations
\begin{equation}\label{jr3}
\begin{array}{rcrcl}
[h_{\theta},\,e_{-i}]&=&t_{i}\,e_{-i},
\qquad [h_{\theta},\,e_{i}]&=&(1-t_{i})\,
e_{i}\quad (t_{i}\in{\mathbb C}),
\\[7pt]
[e_{i},\,e_{-j}]&=&\delta_{ij}\,e_{\theta},
\quad\;[e_{\pm i},\,e_{\pm j}]&=&0~,\quad\;
[e_{\pm i},\,e_{\theta}]=0~,
\end{array}
\end{equation}
provided that
\begin{equation}\label{jr4}
\deg(e_{\theta})=\deg(e_{i})+\deg(e_{-i})
\quad(\mathop{\rm mod}2).
\end{equation}
For the Lie superalgebra $\mathfrak{g}$ the brackets
$[\cdot,\,\cdot]$ always denote the super-commutator:
\begin{equation}\label{jr5}
[x,y]:=xy - (-1)^{\deg(x)\deg(y)}yx
\end{equation}
for any homogeneous elements $x$ and $y$.

Consider the even skew-symmetric two-tensor
\begin{equation}\label{jr6}
r_{\theta,N}^{}(\xi)=\xi\,\Bigl(h_{\theta}\wedge e_{\theta}+
\sum\limits_{i=1}^{N}(-1)^{\deg(e_{i})\,
\deg(e_{-i})}e_{i}\wedge e_{-i}\Bigr)
\end{equation}
where
\begin{equation}\label{jr7}
\deg(\xi)=\deg(e_{\theta}^{})=\deg(e_{i})+\deg(e_{-i})
\quad(\mathop{\rm mod}2),
\end{equation}
and we assume that the operation $"\wedge"$ in (\ref{jr6}) is graded:
\begin{equation}\label{jr8}
e_{i}\wedge e_{-i}:=e_{i}\otimes e_{-i}
-(-1)^{\deg(e_{i})\,\deg(e_{-i})}e_{-i}\otimes e_{i}.
\end{equation}
It is not hard to check that the element (\ref{jr6}) satisfies
the classical Yang-Baxter equation (CYBE),
\begin{equation}\label{jr9}
[r_{\theta,N}^{12}(\xi),r_{\theta,N}^{13}(\xi)+r_{\theta,N}^{23}(\xi)] +
[r_{\theta,N}^{13}(\xi),r_{\theta,N}^{23}(\xi)]=0~,
\end{equation}
and it is called the extended Jordanian $r$-matrix of $N$-order.
Let $N$ be maximal order, i.e. we assume that another elements
$e_{\pm j}\in \mathfrak{n}_{+}$, $j> N$, which satisfy the relations
(\ref{jr3}), do not exist. Such element (\ref{jr6}) is called the
extended Jordanian $r$-matrix of maximal order\cite{T}.

Consider a maximal subalgebra $\mathfrak{b}_+'\in \mathfrak{b}_+$ which
co-commutes with the maximal extended Jordanian $r$-matrix (\ref{jr6}),
$\mathfrak{b}_+':=\mathop{\rm Ker}\delta\in\mathfrak{b_+}$:
\begin{equation}\label{jr10}
\xi\delta(x)=[x\otimes1+1\otimes x,\,r_{\theta,N}^{}(\xi)]=
[\Delta(x),\,r_{\theta,N}^{}(\xi)]=0
\end{equation}
for $\forall x\in\mathfrak{b_+'}$.
Let $r_{\theta_1,N_1}^{}(\xi_1)\in\mathfrak{b}_+'\otimes
\mathfrak{b}_+'$ is also a extended Jordanian $r$-matrix of the form
(\ref{jr5}) with a maximal root $\theta_1\in\mathfrak{h}'$ and
maximal order $N_1$. Then the sum
\begin{equation}\label{jr11}
r_{\theta,N;\,\theta_1,N_1}(\xi,\xi_1^{}):=
r_{\theta,N}^{}(\xi)+r_{\theta_1,N_1}^{}(\xi_{1}^{})
\end{equation}
is also a classical $r$-matrix.

Further, we consider a maximal subalgebra $\mathfrak{b}_+^{\prime\prime}\in
\mathfrak{b}_+'$ which co-commutes with the maximal extended Jordanian
$r$-matrix $r_{\theta_1,N_1}^{}(\xi_{1}^{})$ and we construct a extended
Jordanian $r$-matrix of maximal order, $r_{\theta_{2},N_2}^{}(\xi_{2}^{})$.
Continuing this process as result we obtain a canonical chain of subalgebras
\begin{equation}\label{jr12}
\mathfrak{b}_+\supset\mathfrak{b}_+'\supset\mathfrak{b}_+^{\prime\prime}
\cdots\supset\mathfrak{b}_+^{(k)}
\end{equation}
and the resulting $r$-matrix
\begin{equation}\label{jr13}
r_{\theta,N;\,\ldots;\theta_k,N_k}^{}
\left(\xi,\xi_{1}^{},\cdots,\xi_{k}^{}\right)=
r_{\theta,N}^{}(\xi)+r_{\theta_1,N_1}^{}(\xi_{1}^{})+
\cdots+r_{\theta_k,N_k}^{}(\xi_{k}^{}).
\end{equation}
If the chain (\ref{jr12}) is maximal, i.e. it is
constructed in corresponding with the maximal orders $N, N_1,\dots N_k$,
then the $r$-matrix (\ref{jr13}) is called the maximal classical $r$-matrix
of Jordanian type for the Lie superalgebra $\mathfrak{g}$.

\section{Multiparameter twists of Jordanian type}

The twisting two-tensor $F_{\theta,N}^{}(\xi)$ corresponding to
the $r$-matrix (\ref{jr6}) has the form
\begin{equation}\label{jt1}
F_{\theta,N}(\xi)={\mathcal{F}}_{N}(\xi)F_{\!J}(\sigma_{\theta}^{}),
\end{equation}
where the two-tensor $F_{\!J}^{}$ is the Jordanian twist
and ${\mathcal{F}}_{N}^{}$ is extension of the Jordanian twist
(see\cite{T}). These two-tensors are given by the formulas
\begin{equation}\label{jt2}
F_{\!J}(\sigma_{\theta}^{})=\exp(2h_{\theta}^{}\otimes\sigma_{\theta}^{})~,
\end{equation}
\begin{equation}
\begin{array}{rcl}\label{jt3}
{\mathcal{F}}_{N}(\xi)&=&\Bigl(\prod\limits_{i=1}^{N'}
\exp\bigl(\xi(-1)^{\deg(e_{i})\deg(e_{-i})}e_{i}\otimes e_{-i}\;
e^{-2t_{i}\sigma_{\theta}^{}}\big)\Bigr){\mathcal{F}}_s(\sigma_{\theta}^{})
\\[9pt]
&=&\exp\Big(\xi\sum\limits_{i=1}^{N'} (-1)^{\deg(e_{i})\deg(e_{-i})}
e_{i}\otimes e_{-i}\;e^{-2t_{i}\sigma_{\theta}^{}}\Big)
\,{\mathcal{F}}_s(\sigma_{\theta}^{}),
\end{array}
\end{equation}
where
\begin{equation}\label{jt4}
{\mathcal{F}}_s(\sigma_{\theta}^{})=\Bigl(1-\xi\frac{e_{\theta/2}^{}}
{e^{\sigma_{\theta}}+1}\otimes\frac{e_{\theta/2}^{}}
{e^{\sigma_{\theta}}+1}\Bigr)
\sqrt{\frac{(e^{\sigma_\theta}+1)\otimes(e^{\sigma_\theta}+1)}
{2(e^{\sigma_\theta}\otimes e^{\sigma_\theta}+1)}}~,
\end{equation}
if $\theta/2$ is a root, $e_{\theta/2}^2=e_{\theta}^{}$, $N'=N-1$, and
\begin{equation}\label{jt5}
{\mathcal{F}}_s(\sigma_{\theta}^{})=1~,
\end{equation}
if $\theta/2$ is not any root, $N'=N$. Moreover
\begin{equation}\label{jt6}
\deg(\xi)=\deg(e_{\theta}^{})=\deg(e_{_i})+\deg(e_{-i})
\quad(\mathop{\rm mod}2)~,
\end{equation}
\begin{equation}\label{jt7}
\sigma_{\theta}^{}:=\frac{1}{2}\ln(1+\xi e_{\theta}^{}).
\end{equation}
It should be noted that if the root vector $e_{\theta}^{}$ is odd
then $\sigma_{\theta}^{}=\frac{1}{2}\xi e_{\theta}^{}$.

We can check that the twisting two-tensor (\ref{jt1}) defined by
the formulas (\ref{jt2})--(\ref{jt7}) satisfies the cocycle equation
\begin{equation}\label{jt8}
F^{12}(\Delta\otimes{\rm id})(F)=F^{23}({\rm id}\otimes\Delta)(F)
\end{equation}
and the "unital" normalization condition
\begin{equation}\label{jt9}
(\epsilon \otimes{\rm id})(F)=({\rm id}\otimes\epsilon )(F)=1~.
\end{equation}

The twisted coproduct $\Delta_{\xi}(\,\cdot\,):=
F_{\theta,N}^{}(\xi)\Delta(\,\cdot\,)F_{\theta,N}^{-1}(\xi)$ and
the corresponding antipode $S_{\xi}^{}$ for elements in
(\ref{jr3}) are given by the formulas
\begin{eqnarray}
\Delta_{\xi}(e^{\pm\sigma_{\theta}^{}})&=&
e^{\pm\sigma_{\theta}^{}}\otimes e^{\pm\sigma_{\theta}^{}},\quad
\Delta_{\xi}(e_{\theta/2}^{})\;=\;e_{\theta/2}^{}\otimes1+
e^{\sigma_{\theta}^{}}\otimes e_{\theta/2}^{}~,\quad
\label{jt10}
\\[5pt]
\Delta_{\xi}(h_{\theta}^{})&=&
h_{\theta}^{}\otimes e^{-2\sigma_{\theta}^{}}+1\otimes h_{\theta}^{}+
\frac{\xi}{4}e_{\theta/2}^{}e^{-\sigma_{\theta}^{}}\otimes
e_{\theta/2}^{}e^{-2\sigma_{\theta}^{}}\nonumber
\\[-4pt]
&&-\xi\sum\limits_{i=1}^{N'}(-1)^{\deg e_{i}\deg e_{-i}}
e_{i}\otimes e_{-i}e^{-2(t_{\gamma_i}^{}+1)\sigma_{\theta}^{}},
\label{jt11}
\\[5pt]
&&\;\;\Delta_{\xi}(e_{i})=
e_{i}\otimes e^{-2t_{i}\sigma_{\theta}^{}}+
1\otimes e_{i}, \label{jt12}
\\[5pt]
&&\Delta_{\xi}(e_{-i})=e_{-i}\otimes e^{2t_i\sigma_{\theta}^{}}+
e^{2\sigma_{\theta}^{}}\otimes e_{-i}, \label{jt13}
\end{eqnarray}
\begin{eqnarray}
&&S_{\xi}(e^{\pm\sigma_{\theta}^{}})=e^{\mp\sigma_{\theta}^{}},\qquad
S_{\xi}(e_{\theta/2}^{})=-e_{\theta/2}^{}e^{-\sigma_{\theta}^{}},
\label{jt14}
\\[1pt]
S_{\xi}(h_{\theta}^{})&=&-h_{\theta}^{}e^{2\sigma_{\theta}^{}}+
\frac{1}{4}\bigl(e^{2\sigma_{\theta}^{}}-1\bigr)-
\xi\sum_{i=1}^{N'}(-1)^{\deg(e_{i})\deg(e_{-i})}e_{i}\,e_{-i},
\label{jt15}
\\[1pt]
&&S_{\xi}(e_{_i})=-e_{i}e^{2t_{i}\sigma_{\theta}
},\quad
S_{\xi}(e_{-i})\,=\,-e_{-i}e^{-2(t_{i}+1)\sigma_{\theta}^{}}.
\label{jt16}
\end{eqnarray}
If $\theta/2$ is not any root, the third term in (\ref{jt11}) and
the second term in (\ref{jt15}) should be removed.

The twisted deformation of $U(\mathfrak{g})$ with the new coproduct
$\Delta_{\xi}(\,\cdot\,)$ and the antipode $S_{\xi}^{}$
is  denoted by $U_{\xi}(\mathfrak{g})$.

In order to construct the twist corresponding to the $r$-matrix
(\ref{jr9}) we can not apply the second twist
$F_{\theta_1,N_1}^{}(\xi_1)$ directly in the form
(\ref{jt1})--(\ref{jt4}) to the twisted superalgebra
$U_{\xi}(\mathfrak{g})$ because the deformed coproduct for the
elements of subalgebra $\mathfrak{b}_+'$ can be not trivial, i.e.
\begin{equation}\label{jt17}
\Delta_{\xi}(x)=x\otimes1+1\otimes x + {\rm something}\,,\quad
x\in\mathfrak{b}_+'.
\end{equation}
However, there exists a similarity automorphism $w_{\xi}$ which
trivializes (makes trivial) the twisted coproduct
$\Delta_{\xi}(\,\cdot\,)$ for elements of the subalgebra
$\mathfrak{b}_+'$, i.e.
\begin{equation}\label{jt18}
\Delta_{\xi}(w_{\xi}^{}x w_{\xi}^{-1}):=w_{\xi}^{}x w_{\xi}^{-1}
\otimes1+1\otimes w_{\xi}^{}x w_{\xi}^{-1},\quad
x\in\mathfrak{b}_+'.
\end{equation}
The automorphism $w_{\xi}$ is connected with the Hopf "folding" of
the two-tensor (\ref{jt3}) and it is given by the following formula
(see\cite{T}):
\begin{equation}\label{jt19}
w_{\xi}^{}=\exp\Bigl(\frac{-\xi\sigma_{\theta}}{e^{2\sigma_{\theta}}-1}
\sum_{i=1}^{N'}(-1)^{\deg(e_{i})\deg(e_{i})}e_{i}e_{-i}\Bigr)w_s,
\end{equation}
where $w_s=\exp(\frac{1}{4}\sigma_{\theta})$ if $\theta/2$
is a root, and $w_s=1$ if $\theta/2$ is not any root.

With the help of the automorphism $w_{\xi}$ the total twist chain
corresponding to the $r$-matrix (\ref{jr11}) can be presented as follows
\begin{equation}\label{jt20}
F_{\theta,N;\theta_1,N_1}(\xi,\xi_1)=F_{\theta_{1},N_{1}}(\xi;\xi_1)
F_{\theta,N}^{}(\xi)~,
\end{equation}
where
\begin{equation}\label{jt21}
F_{\theta_1,N_1}(\xi;\xi_1):=(w_{\xi}\otimes w_{\xi})
F_{\theta_1,N_1}(\xi_1)(w_{\xi}^{-1}\otimes w_{\xi}^{-1}).
\end{equation}
Here the two-tensors $F_{\theta,N}^{}(\xi)$ and $F_{\theta_1,N_1}(\xi_1)$
are given by the formulas of type (\ref{jt1})--(\ref{jt5}).

Iterating the formula (\ref{jt21}) we obtain the total twist
corresponding to the $r$-matrix (\ref{jr13}):
\begin{equation}\label{jt22}
\begin{array}{rcl}
&&F_{\theta,N;\theta_1,N_1;\ldots;\theta_k,N_k}(\xi,\xi_1,\ldots,
\xi_k)=F_{\theta_k,N_k}(\xi,\xi_1,\ldots,\xi_{k-1};\xi_k)\cdots
\\[9pt]
&&\qquad\qquad\qquad\qquad\qquad\quad
\times F_{\theta_{2},N_{2}}(\xi,\xi_1;\xi_2)
F_{\theta_{1},N_{1}}(\xi;\xi_1)F_{\theta,N}(\xi)~,
\end{array}
\end{equation}
where $(i=1,\ldots,k)$
\begin{equation}\label{jt23}
\begin{array}{rcl}
&&\!\!F_{\theta_i,N_i}(\xi,\xi_1,\ldots,\xi_{i-1};\xi_i):=
(w_{\xi_{i-1}}\!\!\otimes w_{\xi_{i-1}})\cdots(w_{\xi_{1}}\!\!\otimes
w_{\xi_{1}})(w_{\xi}\!\otimes w_{\xi})
\\[9pt]
&&\qquad\qquad\times\;
F_{\theta_i,N_i}(\xi_i) (w_{\xi}^{-1}\!\otimes w_{\xi}^{-1})
(w_{\xi_{1}}^{-1}\!\otimes w_{\xi_{1}}^{-1})\cdots (w_{\xi_{i-1}}^{-1}
\!\otimes w_{\xi_{i-1}}^{-1})\,.
\end{array}
\end{equation}

Now we consider specifically the multiparameter twists for the classical
superalgebras $\mathfrak{gl}(m|n)$ and $\mathfrak{osp}(1|2n)$.

\section{Quantum deformation of Jordanian type for 
$\mathfrak{gl}(m|n)$}

Let $e_{ij}$ ($i,j=1,2,\ldots,m+n$) be standard
$(n+m)\times(n+m)$-matrices, where $(e_{ij})_{kl}=\delta_{ik}\delta_{jl}$.
For such matrices we define a supercommutator as follows
\begin{equation}\label{gl1}
[e_{ij},\,e_{kl}]:=e_{ij}e_{kl}-(-1)^{\deg(e_{ij})\deg(e_{kl})}e_{kl}e_{ij},
\end{equation}
where $\deg(e_{ij})=0$ for $i,j\leq n$ or $i,j> n$, and $\deg(e_{ij})=0$
in another cases. It is easy to check that
\begin{equation}\label{gl2}
[e_{ij},\,e_{kl}]=\delta_{jk}e_{il}-
(-1)^{\deg(e_{ij})\deg(e_{kl})}\delta_{il}e_{kj}.
\end{equation}
The elements $e_{ij}$ ($i,j=1,2,\ldots,N:=m+n$) with the relations
(\ref{gl2}) are generated the Lie superalgebra $\mathfrak{gl}(m|n)$.

The maximal $r$-matrix of Jordanian type for the Lie superalgebra
$\mathfrak{gl}(m|n)$ has the form{\cite{T}
\begin{equation}\label{gl3}
r_{1,\ldots,[N/2]}^{}
\bigl(\xi_{1}^{},\ldots,\xi_{[N/2]}^{}\bigr)=r_{1}^{}(\xi_{1}^{})
+\cdots+r_{[N/2]}^{}(\xi_{[N/2]}^{})\,,
\end{equation}
where ($i=1,2,\ldots, [N/2]$)
\begin{equation}\label{gl4}
\begin{array}{rcl}
r_{i}(\xi_{i})&=&{\displaystyle\xi_{i}\Bigl(\frac{1}{2}}
(e_{ii}^{}-e_{N+1-i,N+1-i}^{})\wedge e_{i,N+1-i}^{}+
\\[9pt]
&&\quad{\displaystyle+\sum_{k=i+1}^{N-i}}
(-1)^{\deg(e_{ik}^{})\deg(e_{kN+1-i}^{})}e_{ik}\wedge e_{kN+1-i}\Bigr)~.
\end{array}
\end{equation}
Consider the first twist corresponding to the $r$-matrix
$r_{1}(\xi_{1})$
\begin{equation}\label{gl5}
F_{1,N-2}(\xi_1)\;=\;{\mathcal{F}}_{N-2}(\xi_1)F_{\!J}(\sigma_{1}^{}),
\end{equation}
where
\begin{equation}\label{gl6}
F_{\!J}^{}(\sigma_{1}^{})= e^{(e_{11}^{}-e_{NN}^{})\otimes\sigma_{1}^{}},
\end{equation}
\begin{equation} \label{gl7}
{\mathcal{F}}_{N-2}(\xi_1)=\exp\Big(\xi_1\sum\limits_{k=2}^{N-1}
(-1)^{\deg(e_{1k})\deg(e_{kN})}
e_{1k}\otimes e_{kN}\;e^{-2\sigma_{1}^{}}\Big),
\end{equation}
\begin{equation}\label{gl8}
\sigma_{1}^{}:=\frac{1}{2}\ln(1+\xi_1 e_{1N}^{}).
\end{equation}
The corresponding automorphism $w_{\xi}$ is connected with the Hopf
"folding" of the two-tensor (\ref{gl7}) and is given
as follows
\begin{equation}\label{gl9}
w_{\xi_1^{}}=\exp\Bigl(\frac{-\xi_1\sigma_{1}^{}}{e^{2\sigma_{1}^{}}-1}
\sum_{k=2}^{N-1}(-1)^{\deg(e_{1k}^{})\deg(e_{kN}^{})}e_{1k}^{}e_{kN}^{}\Bigr).
\end{equation}
It is easy to see that
\begin{equation}\label{gl10}
w_{\xi_1}e_{ij}w_{\xi_1}^{-1}=e_{ij}
\end{equation}
for all $i,\,j$ satisfying the condition $2\leq i,\,j\leq N-2$, therefore
(see the formula (\ref{jt18})) deformed coproducts
$\Delta_{\xi_1}(\,\cdot\,):=F_{1,N-2}^{}(\xi_1)\Delta(\,\cdot\,)
F_{1,N-2}^{-1}(\xi_1)$ for these elements are trivial:
\begin{equation}\label{gl11}
\Delta_{\xi_1}(e_{ij})=e_{ij}\otimes1+1\otimes e_{ij},\quad
2\leq i,\,j\leq N\!-\!1.
\end{equation}
This means that the automorphism $w_{\xi_1}$ in the formula
(\ref{jt21}) for the case $\mathfrak{gl}(m|n)$ acts trivially and
therefore the total twist corresponding to the $r$-matrix (\ref{gl3})
is given as follows
\begin{equation}\label{gl12}
\begin{array}{rcl}
&&F_{1,N-2;2,N-4;\ldots;k,N-2k}(\xi_1,\xi_2,\ldots,
\xi_k)=F_{k,N-2k}(\xi_k)\cdots
\\[9pt]
&&\qquad\qquad\qquad\qquad\qquad\qquad\qquad\quad
\times F_{{2},N-4}(\xi_2)F_{{1},N-2}(\xi_1)~,
\end{array}
\end{equation}
where $(i=1,\ldots,[N/2])$
\begin{equation} \label{gl13}
\begin{array}{rcl}
F_{i,N-2i}(\xi_i)&=&\exp\Big(\xi_i\!\!\sum\limits_{k=i+1}^{N-i}
\!(-1)^{\deg(e_{ik})\deg(e_{kN-2i})}e_{ik}\otimes e_{kN-2i}\,
e^{-2\sigma_{i}^{}}\Big)
\\[12pt]
&&\times\exp\big((e_{ii}^{}-e_{N-2i,N-2i}^{})\otimes
\sigma_{i}^{}\big)\,.
\end{array}
\end{equation}

\section{Quantum deformation of Jordanian type for 
$\mathfrak{osp}(1|2n)$}

In order to obtain compact formulas describing the commutation
relations for generators of the orthosymplectic superalgebra
$C(n)\simeq\mathfrak{osp}(1|2n)$} we use embedding of this superalgebra
in the general linear superalgebra $\mathfrak{gl}(1|2n)$.
Let $a_{ij}$ ($i,j=0,\pm1,\pm2, \ldots,\pm n$) be a standard
basis of the superalgebra $\mathfrak{gl}(1|2n)$ (see the previous
Section 4) with the standard supercommutation relations
\begin{equation}\label{cw1}
[a_{ij},\,a_{kl}]=\delta_{jk}a_{il}-(-1)^{\deg(e_{ij})\deg(e_{kl})}
\delta_{il}a_{kj},
\end{equation}
where $\deg(e_{ij})=1$ when one index $i$ or $j$ is equal to 0
and another takes any value $\pm1,\ldots,\pm n$;
 $\deg(e_{ij})=0$ in the remaining  cases. The superalgebra
$\mathfrak{osp}(1|2n)$ is embedded in $\mathfrak{gl}(1|2n)$ as a
linear envelope of the following generators:
\\(i) the even (boson) generators spanning  the symplectic algebra
$\mathfrak{sp}(2n)$:
\begin{equation}\label{cw2}
e_{ij}:=a_{i-j}+\sign(ij)\,a_{j-i}=
\sign(ij)\,e_{ji}^{}\quad(i,j=\pm1,\pm2,\ldots,\pm n)~;
\end{equation}
(ii) the odd (fermion) generators extending $\mathfrak{sp}(2n)$
 to $\mathfrak{osp}(1|2n)$:
\begin{equation}\label{cw3}
e_{0i}:=a_{0-i}+\sign(i)\,a_{i0}=
\sign(i)\,e_{i0}\quad(i=\pm1,\pm2,\ldots,\pm n)~.
\end{equation}
We also set $e_{00}^{}=0$ and introduce the sign function: $\sign
x=1$ if a real number $x\geq0$ and $\sign x=-1$ if $x<0$. One can
check  that the elements (\ref{cw2}) and (\ref{cw3}) satisfy the
following relations:
\begin{eqnarray}\label{cw4}
[e_{ij},\,e_{kl}]&=&\delta_{j-k}e_{il}+\delta_{j-l}
\sign(kl)\,e_{ik}-\delta_{i-l}e_{kj}-\delta_{i-k}\sign(kl)\,e_{lj},
\\[3pt]
[e_{ij},\,e_{0k}]&=&\delta_{j-k}\sign(k)\,e_{i0}-
\delta_{i-k}e_{0j},\label{cw5}
\\[3pt]
\{e_{0i},e_{0k}\}&=&\sign(i)\,e_{ik}
\label{cw6}
\end{eqnarray}
for all $i,j,k,l=\pm1,\pm2,\ldots,\pm n$, where the bracket
$\{\cdot,\cdot\}$ means anti-commutator.

The elements $e_{ij}^{}$ $(i,j=0,\pm1,\pm2)$ are not linearly
independent (we have for example, $e_{1-2}^{}=-e_{-21}^{})$ and we
can choose from them the Cartan-Weyl basis as follows
\begin{eqnarray}\label{cw7}
&&\!\!\!\!\!\!\!\!
{\rm rising\; generators}:
e_{i\pm j},e_{kk},e_{0k}\;(1\leq i<j\leq n,\;1\leq k\leq n);
\\[3pt]\label{cw8}
&&\!\!\!\!\!\!\!\!
{\rm lowering\; generators}:
e_{\pm j-i},e_{-k-k},e_{-k0}\;(1\leq i<j\leq n,\;1\leq k\leq n);
\\[3pt]\label{cw9}
&&\!\!\!\!\!\!\!\!
{\rm Cartan\; generators}:h_{i}:=
e_{k-k}\;\; (1\leq k\leq n).
\end{eqnarray}

Maximal classical r-matrix of Jordanian type for the Lie superalgebra
$\mathfrak{osp}(1|2n)$ has the form{\cite{T}
\begin{equation}\label{cw10}
r_{1,2,\ldots,n}^{}\left(\xi_{1},\xi_{2}^{},\cdots,\xi_{n}^{}\right)=
r_{1}^{}(\xi_{1}^{})+r_{2}^{}(\xi_{2}^{})+\cdots+r_{n}^{}(\xi_{n}^{})~.
\end{equation}
where
\begin{equation}\label{cw11}
\begin{array}{rcl}
r_{i}^{}\left(\xi_{i}^{}\right)&:=&\xi_{i}^{}\Bigl
(\displaystyle{\frac{1}{2}}e_{i-i}^{}\wedge e_{ii}-
2e_{0i}\otimes e_{0i}+\sum_{k=i+1}^{n}e_{i-k}\wedge e_{ik}\Bigr)~,
\end{array}
\end{equation}
The total twist corresponding to the $r$-matrix (\ref{cw10}) is given
as follows
\begin{equation}\label{cw12}
\begin{array}{rcl}
F_{1,n;2,n-1;\ldots;n,1}(\xi,\xi_1,\ldots,\xi_n)&=&
F_{n,1}(\xi_1,\xi_2,\ldots,\xi_{n-1};\xi_n)\cdots
\\[9pt]
&&\times F_{2,n-1}(\xi_1;\xi_2)F_{1,n}(\xi_1)~.
\end{array}
\end{equation}
Here $(i=1,\ldots,k)$
\begin{equation}\label{cw13}
\begin{array}{rcl}
&&\!\!F_{i,n+1-i}(\xi_1,\ldots,\xi_{i-1};\xi_i):=
(w_{\xi_{i-1}}\!\!\otimes w_{\xi_{i-1}})\cdots
w_{\xi_{2}})
(w_{\xi_1}\!\otimes w_{\xi_1})
\\[7pt]
&&\qquad\qquad\qquad\times\;
F_{i,n+1-i}(\xi_i) (w_{\xi_1}^{-1}\!\otimes w_{\xi_1}^{-1})
\cdots (w_{\xi_{i-1}}^{-1}\!\otimes w_{\xi_{i-1}}^{-1})\,.
\end{array}
\end{equation}
\begin{equation} \label{cw14}
F_{i,n+1-i}(\xi_i)=\exp\Big(\xi_i\!\!\sum\limits_{k=i+1}^{n+1-i}
e_{i-k}\otimes e_{k,n+1-i}\,e^{-2\sigma_{i}^{}}\Big)
{\mathcal{F}}_s(\sigma_{i}^{})\,e^{e_{i-i}\otimes\sigma_{i}^{}},
\end{equation}
where $\mathcal{F}_s(\sigma_{i}^{})$ is defined by the formula
(\ref{jt4}), and
\begin{equation}\label{cw15}
w_{\xi_i^{}}=\exp\Bigl(\frac{-\xi_i\sigma_{i}^{}}{e^{2\sigma_{i}^{}}-1}
\sum_{k=i+1}^{n+1-i}e_{i-k}^{}e_{k,n+1-i}^{}\Bigr),\quad
\sigma_{i}^{}:=\frac{1}{2}\ln(1+\xi_1 e_{i-i}^{}).
\end{equation}

\section{Light-cone $\kappa$-deformation of the super-Poincar\'{e}\\
algebra ${\mathcal P}(3,1|1)$}

The Poincar\'{e} algebra ${\mathcal{P}}(3,1)$ of the 4-dimensional
space-time is generated by 10 elements, $M_j$,  $N_j$, $P_j$, $P_0$
$(j=1,2,3)$ with the standard commutation relations:
\begin{equation}
\begin{array}{rcccccl}\label{lc1}
[M_j,\,M_k]&=&\!i\epsilon_{jkl}\,M_l,\quad
[M_j,\,N_k]&=&\!i\epsilon_{jkl}\,N_l,\quad
[N_j,\,N_k]&=&\!-i\epsilon_{jkl}\,M_l,
\\[5pt]
[M_j,\,P_k]&=&\! i\epsilon_{jkl}\,P_l,\quad\,
[M_j,\,P_0]&=&0,\phantom{aaaaaaaaaaaaa}&&
\\[5pt]
[N_j,\,P_k]&=&\!-i\delta_{jk}\,P_0,\quad\;
[N_j,\,P_0]&=&\!-iP_j^{},\qquad
[P_\mu,\,P_\nu]&=&0~.
\end{array}
\end{equation}
The super-Poincar\'{e} algebra ${\mathcal P}(3,1|1)$ is generated
by the algebra $\mathcal{P}(3,1)$ and four real supercharges
$Q_\alpha^{}$ $(\alpha=\pm1,\pm2)$) with the commutation relations
\begin{equation}
\begin{array}{rcl}\label{lc2}
[M_j^{},Q^{(\pm)}_{\alpha}]&=&-\frac{i}{2}(\sigma_j^{})_
{\alpha\beta}^{}\,Q^{(\pm)}_{\beta},\quad
\\[5pt]
[N_j^{},Q^{(\pm)}_{\alpha}]&=&\mp\frac{i}{2}(\sigma_j^{})_
{\alpha\beta}^{}\,Q^{(\pm)}_{\beta},\quad
[P_\mu^{},Q^{(\pm)}_{\alpha}]\,=\,0,
\end{array}
\end{equation}
and moreover
\begin{equation}
\label{lc3}
\{Q^{(\pm)}_{\alpha},\,Q^{(\pm)}_{\beta}\}=0,\quad
\{Q^{(+)}_{\alpha},\,Q^{(-)}_{\beta}]=2\bigl(
\delta_{\alpha\beta}^{}\,P_0^{}-
(\sigma_j^{})_{\alpha\beta}^{}\,P_j^{}\bigr),
\end{equation}
where we use the denotations $Q^{(\pm)}_1:=Q_1^{}\pm iQ_2^{}$,
$Q^{(\pm)}_2:=Q_{-1}^{}\pm iQ_{-2}^{}$, and $\sigma_j^{}$ $(j=1,2,3)$
are $2\times2$ $\sigma$-matrices. It should be noted that
the spinor ${\bf Q}^{(+)}:=(Q^{(+)}_{1},Q^{(+)}_{2})$ transformes
as the left-regular representation and the spinor ${\bf
Q}^{(-)}:=(Q^{(-)}_{1},Q^{(-)}_{2})$ provides the right-regular
one with respect to ${\mathcal{P}}(3,1)$.

Using the commutation relations (\ref{lc1}) and (\ref{lc2}),
(\ref{lc3}) it is easy to check that the elements $iN_3$,
$P_+:=P_0+P_3$, $P_1$, $i(N_1+M_2)$, $P_2$, $i(N_2-M_1)$, and $Q_{\alpha}$
($\alpha=1,2$) satisfy the relations (\ref{jr1})--(\ref{jr4}), namely,
$\{h_{\gamma_0^{}},e_{\gamma_0^{}}\}\rightarrow\{iN_3,P_+\}$,
$\{e_1,e_{-1}\}\rightarrow\{P_1,i(N_1+M_2)\}$,
$\{e_2,e_{-2}\}\rightarrow\{P_2,i(N_2-M_1)\}$,
$e_{\pm 3}\rightarrow Q_1$, $e_{\pm 4}\rightarrow Q_2$.
Therefore the  two-tensor
\begin{equation}\begin{array}{rcl}\label{lc4}
r&=&\frac{1}{\kappa}\Bigl(P_1^{}\wedge (N_1^{}+M_2^{})
+P_2{}\wedge (N_2^{}-M_1^{})+P_+\wedge N_3^{}+
\\[5pt]
&&\quad\;\;+2\bigl(Q_1^{}\wedge Q_1^{}+Q_2^{}\wedge Q_2^{}\bigr)\Bigr)~,
\end{array}
\end{equation}
is a classical $r$-matrix of Jordanian type. It is called the
classical $r$-matrix for light-cone $\kappa$-deformation of $D=4$
super-Poincar\'{e}. Specializing the general formula (\ref{jt3})
to our case ${\mathcal P}(3,1|1)$ we immediately obtain the twisting
two-tensor corresponding to this $r$-matrix
\begin{eqnarray}\label{lc5}
F_{\kappa}^{}({\mathcal P}(3,1|1)&:=&
\mathfrak{F}_{\kappa}^{}(Q_2^{})\mathfrak{F}_{\kappa}^{}(Q_1^{})
F_{\kappa}^{}({\mathcal P}(3,1))~,
\end{eqnarray}
where $F_{\kappa}^{}({\mathcal P}(3,1))$ is the twisting two-tensor of
the light-cone $\kappa$-deformation of the Poinca\-r\'{e} algebra
${\mathcal P}(3,1)$
\begin{eqnarray}\label{lc6}
F_{\kappa}^{}({\mathcal P}(3,1))\!\!&:=\!\!&e^{\frac{i}{\kappa}\,
P_1^{}\otimes (N_1^{}+M_2^{})e^{-2\sigma_{\!+}^{}}}\,e^{\frac{i}
{\kappa}\,P_2^{}\otimes (N_2^{}-M_1^{})e^{-2\sigma_{\!+}^{}}}
\,e^{2iN_3\otimes\sigma_{\!+}^{}}
\end{eqnarray}
and the super-factors $\mathfrak{F}_{\kappa}^{}(Q_\alpha)$
($\alpha=1,2$) are given by the formula
\begin{equation}\label{lc7}
\mathfrak{F}_\kappa^{}(Q_\alpha^{})=\sqrt{\frac{(1+
e^{\sigma_{\!+}^{}})\otimes(1+e^{\sigma_{\!+}^{}})}
{2(1+e^{\sigma_{\!+}^{}}\!\otimes e^{\sigma_{+}^{}})}}
\biggl(1+\frac{2}{\kappa}\,\frac{Q_{\alpha}^{}}
{1+e^{\sigma_{\!+}^{}}}\otimes\frac{Q_{\alpha}^{}}
{1+e^{\sigma_{\!+}^{}}}\biggr)~,
\end{equation}
\begin{equation}\label{lc8}
\sigma_{\!+}^{}:=\frac{1}{2}\ln\Bigl(1+\frac{1}{\kappa}P_{+}^{}\Bigr)~.
\end{equation}
The formulas (\ref{lc5})--(\ref{lc8}) were obtained by a suitable
contraction of the quantum deformation of Jordanian type\cite{BLT}.


\end{document}